%% file: smoothpoly.tex
\documentstyle[12pt, fullpage]{amsart}
\begin{document}
\thispagestyle{empty}
\input{universal}
\createtitle{Lower Bounds for the Number of Smooth Values
of~a~Polynomial}{11N32}
%\doublespace

\abstract

We investigate the problem of showing that the values of a given
polynomial are smooth (i.e., have no large prime factors) a positive
proportion of the time. Although some results exist that bound the
number of smooth values of a polynomial from above, a corresponding
lower bound of the correct order of magnitude has hitherto been
established only in a few special cases. The purpose of this paper is
to provide such a lower bound for an arbitrary polynomial. Various
generalizations to subsets of the set of values taken by a polynomial
are also obtained.

\endabstract

\def\V{{\cal V}}\def\S{{\cal S}}\def\P{{\cal P}}\def\A{{\cal A}}
\def\lmp#1{\lambda_{#1}}
\def\3#1{\#\{#1\}}

\def\sigmar#1{\sigma_{#1}}
\newsymbol\lsim 132E % see The Joy of TeX (2nd ed.), Appendix G
\def\ord{\mathop{\rm ord}}
\def\sst{\tau^*}
\def\sh{L}
\def\osm#1{$O(#1)$-smooth}
\def\hhh{h_2}\def\hht{h_1}
\def\mps{\mod{p^r}}\def\mpsi{\mod{p^{r+1}}}

\def\gov#1{{g(#1)\over{#1}}}
\def\glov#1{{g(#1)\log{#1}\over{#1}}}
\def\agov#1{{\abs{g(#1)}\over{#1}}}
\def\aglov#1{{\abs{g(#1)}\log{#1}\over{#1}}}
\def\gqov#1{{g_q(#1)\over{#1}}}

\section{Introduction}\label{introsec}\noindent\input{intro}
\section{Outline of the Approach and Notation}\label{grungysec}\noindent\input{grungy}

\section{Multiplicative Functions Associated to a Polynomial}\label{multsec}\noindent\input{mult}

\section{Sums of Multiplicative Functions}\label{mult2sec}\noindent\input{mult2}

\section{Application of the Upper Bound Sieve}
 \Mmake{sieve}\input{Sd}\input{sieve}\input{cancel}\input{apply}
\section{Smooth Values on Prime Arguments}\label{primesec}\noindent\input{prime}
\section{An Elementary Lower Bound}\label{sillysec}\noindent\input{silly}

\bibliography{smoothpoly}
\bibliographystyle{../amsplain}
\end{document}

%% file: universal.tex
\newcommand{\createtitle}[2]{\title{#1}\author{Greg Martin}\address{Department of Mathematics\\University of Toronto\\Canada M5S 3G3}\email{gerg@@math.toronto.edu}\subjclass{#2}\maketitle}
\newcommand{\make}[1]{\label{#1sec}\noindent\input{#1}}
\newcommand{\Mmake}[1]{\label{#1sec}\noindent}

%%% use for AMSTeX
% \input{/opt/texmf/tex/ams/amstex}
%%%

\newtheorem{theorem}{Theorem}
\newtheorem{lemma}[theorem]{Lemma}
\newtheorem{corollary}{Corollary}[theorem]
\newtheorem{proposition}[theorem]{Proposition}

\newenvironment{pflike}[1]{\noindent{\bf #1}}{\vskip10pt} % NO \QED
\newenvironment{proof}{\begin{pflike}{Proof:}}{\qed\end{pflike}}

\def\(#1\)#2{
 \def\helpera{\ifcase#2(\or\big(\or\Big(\or\bigg(\or\Bigg(\else(\fi}
 \def\helperb{\ifcase#2)\or\big)\or\Big)\or\bigg)\or\Bigg)\else)\fi}
 \helpera{#1}\helperb}

\newcommand{\2}[1]{\ifmmode{\cal#1}\else$\cal#1$\fi}
\newcommand{\3}[1]{\#\{#1\}}
\newcommand{\abs}[1]{\left|#1\right|} % take out
\newcommand{\floor}[1]{\lfloor#1\rfloor}
\newcommand{\bfloor}[1]{\big\lfloor#1\big\rfloor}
\newcommand{\bbfloor}[1]{\bigg\lfloor#1\bigg\rfloor}
\newcommand{\ceil}[1]{\lceil#1\rceil}
\newcommand{\bceil}[1]{\big\lceil#1\big\rceil}
\newcommand{\bbceil}[1]{\bigg\lceil#1\bigg\rceil}

\renewcommand{\mod}[1]{{\ifmmode\text{\rm\ (mod~$#1$)}
 \else\discretionary{}{}{\hbox{ }}\rm(mod~$#1$)\fi}}
\newcommand{\ep}{\varepsilon}
\renewcommand{\implies}{\Rightarrow}
\newcommand{\rmif}{{\rm if\ }}

\newcommand{\half}{{\mathchoice{\textstyle\frac12}{1/2}{1/2}{1/2}}}
\newsymbol\dnd 232D % see The Joy of TeX (2nd ed.), Appendix G
\newcommand{\exdiv}{\mathrel{\mid\mid}}

\renewcommand{\lg}[1]{\mathop{\log_{#1}}}
\def\lgs#1^#2{\mathop{\log_{#1}^{#2}}}
\newcommand{\li}{\mathop{\rm li}}

\newcommand{\doublespace}{%\parskip 0.2cm
  \baselineskip=24pt}
\newcommand{\spaceandahalf}{\parskip6pt\baselineskip=18pt}

\newcommand{\scroll}[1]{\scrollmode#1\errorstopmode}
\newcommand{\comment}[1]{}

\vfuzz=2pt % prohibits overfull \vbox messages due to page headers
	   % and footers

%% file: intro.tex
Our knowledge of the multiplicative properties of the values taken by
a polynomial with integer coefficients (or, more generally,
an integer-valued polynomial) is quite limited. For instance, it is
conjectured that if $h(t)$ is a polynomial that is not identically
zero modulo any prime, then the irreducible factors of $h$ will
simultaneously take prime values for infinitely many values of $n$; in
fact, there is a conjectured asymptotic formula (see Bateman--Horn
\cite{BatHor}) for the number of positive integers $n\le x$ for which
this occurs. Dirichlet's theorem on primes in arithmetic progressions
verifies this conjecture when $h$ is a linear polynomial, but when $h$
has degree at least 2, these conjectures are still unresolved; it is
unknown, for instance, whether there are infinitely many primes of the
form $n^2+1$, or whether there are infinitely many primes $p$ such
that $p+2$ is also prime (the twin primes conjecture).

Another multiplicative property of integers is {\it smoothness\/}: an
integer is $y$-smooth if none of its prime factors exceed $y$. Since
an integer $n$ is prime if and only if all of its prime factors exceed
$n^{1/2}$, smoothness is in some sense the complementary property to
being prime. If we define $\Psi(x,y)$ to be the number of $y$-smooth
positive integers not exceeding $x$, then it is well-known that
$\Psi(x,x^{1/u})$ is asymptotic to $\rho(u)x$ for fixed $u$ or for
$u$ growing not too quickly with $x$, where $\rho$ is the solution of
a particular differential-difference equation. In particular, for a
fixed real number $0<\alpha<1$, the $x^\alpha$-smooth integers
comprise a positive proportion of the integers up to $x$.

When $h$ is a polynomial of degree 1, we again have an asymptotic
formula for the number of integers $n\le x$ for which $h(n)$ is
$x^{1/u}$-smooth, which for fixed $u$ and $h$ was first established in
the work of Buchstab \cite{Buc:OtNiaAP} on smooth numbers in
arithmetic progressions (later work has provided results having some
uniformity in the coefficients of the linear polynomial; see
Hildebrand--Tenenbaum \cite[Section 6]{HilTen:IwLPF} for a discussion
of such results). Our qualitative understanding of the smooth values
of a fixed polynomial $h$ of degree $g\ge2$ is somewhat better than
that of its prime values. Schinzel shows \cite[Theorem 13]{Sch:OTToG}
that
\begin{equation}
\hbox{\it there are infinitely many integers $n$ for which $h(n)$ is
$n^{g-1-\delta(g)}$-smooth,}
\label{schin}
\end{equation}
where $\delta(g)$ is a certain real number which satisfies
$0<\delta(g)<1$ and $\delta(g)=2/g+O(g^{-2})$ for large $g$; this is a
nontrivial result because $h(n)$ has order of magnitude $n^g$. He also
shows \cite[Theorem 15]{Sch:OTToG} that if $h$ has the special form
$h(t)=At^g+B$ for some nonzero integers $A$, $B$, and $g\ge2$, then
for any positive real number $\alpha$ there are infinitely many
integers $n$ for which $h(n)$ is $n^\alpha$-smooth. The same
conclusion also holds, by work of Balog and Wooley \cite{BalWoo:OSoCI}
(extending a result of Eggleton and Selfridge \cite{EggSel:CIwNLPF}),
when $h(t)=\prod_{1\le i\le g}(A_it+B_i)$ is the product of linear
polynomials with integer coefficients.

Unfortunately, the proofs of these results do not give very strong
estimates for how many smooth values are taken by $h$. If we define
the counting function of the $y$-smooth values of~$h$,
\begin{equation*}
\Psi(h;x,y)=\#\{1\le n\le x:p\mid h(n)\implies p\le y\}
\end{equation*}
(where $p$ generically denotes a prime), then presumably, for any
fixed polynomial $h$ and positive real number $\alpha$, we should have
$\Psi(h;x,x^\alpha)\sim c(h,\alpha)x$ for some positive constant
$c(h,\alpha)$. However, the arguments of Schinzel and Balog--Wooley
imply only lower bounds of the form $\Psi(h;x,y)\gg x^\beta$ for
rather small values of $\beta$.

When $h$ is a linear polynomial, Buchstab's work referred to above
gives the asymptotic formula $\Psi(h;x,x^\alpha)\sim
\rho(\alpha^{-1})x$, directly extending the formula for
$\Psi(x,x^{1/u})$ mentioned earlier. There are a few results of
Hmyrova \cite{Hmyrova, HmyrovaII} and Timofeev \cite{Tim:PwSPD} that
give upper bounds for $\Psi(h;x,x^\alpha)$ for polynomials of
arbitrary degree; however, there has been very little progress towards
establishing a lower bound of the presumed order of magnitude,
\begin{equation}
\Psi(h,x,x^\alpha)\gg_{h,\alpha} x,  \label{presumed}
\end{equation}
for polynomials of degree at least 2. It is known that the lower bound
(\ref{presumed}) holds for any $\alpha>0$ when $h$ has the form
$h(t)=At(Bt+C)$, where $A$, $B$, and $C$ are integers with $AB>0$, by
work of Balog and Ruzsa \cite{BalRuz:OaAPoSS} (generalizing a result
of Hildebrand \cite{Hil:OaCoB}). It also holds for
$\alpha>e^{-1/(g-1)}$ when $h(t)=(t+1)(t+2)\dots(t+g)$ for some
$g\ge2$, by work of Hildebrand \cite{Hil:OISCSoCI}. The only result
along these lines for irreducible polynomials is due to Dartyge
\cite{Dar:EdlF}, who shows that (\ref{presumed}) holds for
$\alpha>149/179$ when $h(t)=t^2+1$.

We are able to establish a lower bound of the form (\ref{presumed})
for an arbitrary polynomial, as indicated in the following theorem.

\begin{theorem}
Let $h(t)$ be an integer-valued polynomial (not identically zero), and
let $g$ be the largest of the degrees of the irreducible factors of
$h$. Let $k$ be the number of distinct irreducible factors of $h$ of
degree $g$, and let $\delta$ be any positive real number less than
$(2k+1)^{-1}$. Then when $x$ is sufficiently large, we have
\begin{equation}
\Psi(h;x,x^{g-\delta})\gg_{h,\delta}x.  \label{mainbd}
\end{equation}
In particular, if $h$ is irreducible, then the lower bound
{\rm(\ref{mainbd})} holds for any $0<\delta<1/3$.
\label{mainthm}
\end{theorem}
\noindent By the definition of $g$, the values $h(n)$ with $n\le x$
are trivially $O_h(x^g)$-smooth; Theorem \ref{mainthm} asserts that a
positive proportion of the values $h(n)$ with $n\le x$ are
$x^{g-\delta}$-smooth. One feature of this result is that the amount
$x^\delta$ that we are able to save from the trivial smoothness
parameter does not depend on the polynomial $h$, but only on the
degrees of its irreducible factors.

Our methods can be extended to show the abundance of smooth values
$h(n)$ with $n$ restricted to various sets. Our goal is to obtain a
lower bound of the correct order of magnitude for the number of
$y$-smooth values $h(n)$, for some non-trivial value of $y$; it turns
out that we can do this with $n$ restricted in a wide variety of
ways. For the purposes of illustration, we state the following
theorems.

\begin{theorem}
Let $h(t)$, $g$, $k$, and $\delta$ be as in Theorem \ref{mainthm}.
For real numbers $x\ge\sh\ge2$, define
\begin{equation*}
\Psi(h;x,\sh,y) = \Psi(h;x,y) - \Psi(h;x-\sh,y).
\end{equation*}
Then when $x$ is sufficiently large, we have
\begin{equation}
\Psi(h;x,\sh,x^g\sh^{-\delta}) \gg_{h,\delta} \sh.  \label{shortbd}
\end{equation}
In particular, if $h$ is irreducible, then the lower bound
{\rm(\ref{shortbd})} holds for any $0<\delta<1/3$.
\label{shortthm}
\end{theorem}
\noindent Thus a positive proportion of the values taken a polynomial
on a short interval of length $\sh$ are nontrivially smooth by a
fractional power of $\sh$.

\begin{theorem}
Let $h(t)$, $g$, and $k$ be as in Theorem \ref{mainthm}. Let $\A$ be
any set of integers whose density $\eta$ exists and is positive, and
let $\delta$ be any positive real number less than $\eta/(2k+\eta)$.
Define
\begin{equation*}
\Psi_\A(h;x,y)=\#\{1\le n\le x,\, n\in\A:p\mid h(n)\implies p\le y\}.
\end{equation*}
Then when $x$ is sufficiently large, we have
\begin{equation*}
\Psi_\A(h;x,x^{g-\delta}) \gg_{h,\eta,\delta} x.
\end{equation*}
\label{scantthm}
\end{theorem}
\noindent For example, if $h$ is an irreducible polynomial of degree
$g$, then a positive proportion of the values that $h$ takes on
squarefree integers are $x^{g-\delta}$-smooth for any
$\delta<3/(\pi^2+3)=0.2331\dots$. A suitably modified theorem can be
established for sets of integers whose densities do not exist, one
consequence of which is the following: if $\A$ is a set of integers
such that there is never an abundance of values $h(a)$, $a\in\A$, that
are nontrivially smooth by a power of $x$ (i.e., if $\lim_{x\to\infty}
\Psi_\A(h;x,x^{g-\ep})/x =0$ for every $\ep>0$), then $\A$ must have
density 0.

It is worth noting that the proofs of Theorems
\ref{mainthm}--\ref{scantthm} can be extended to the case where $h$ is
a polynomial in more than one variable. In fact, one can obtain
stronger results, in terms of the admissible ranges of the smoothness
parameter, by a more sophisticated treatment of the error terms
arising in the application of the sieve in Section \ref{sievesec}. For
instance, the lower bound (\ref{mainbd}) holds for $\delta$ as large
as $1/2+o(1)$ as the number of variables increases, at least under
some hypothesis controlling the singularities of the polynomial. We do
not discuss the details herein.

The values that a polynomial takes on prime arguments form a natural
arithmetic set, and the question of whether such a set contains
infinitely many prime numbers is an important motivating problem of
sieve theory. For example, when $h(t)=t+2$, this question is precisely
the twin primes conjecture. Analogously, we can ask whether such a set
contains many smooth numbers; the following theorem demonstrates that
it does.

\begin{theorem}
Let $h(t)$, $g$, and $k$ be as in Theorem \ref{mainthm}, and let
$\delta$ be a positive real number less than $(4k+2)^{-1}$. Define
\begin{equation*}
\Phi(h;x,y)=\#\{1\le q\le x,\, q\hbox{ \rm prime} : p\mid h(q)\implies
p\le y\}.
\end{equation*}
Then when $x$ is sufficiently large, we have
\begin{equation}
\Phi(h;x,x^{g-\delta})\gg_{h,\delta}x/\log x.  \label{primebd}
\end{equation}
In particular, if $h$ is irreducible, then the lower bound
{\rm(\ref{primebd})} holds for any $0<\delta<1/6$.
\label{primethm}
\end{theorem}
\noindent Thus a positive proportion of the values a polynomial takes
on primes are nontrivially smooth by a power of $x$. The
aforementioned work of Hmyrova contains upper bounds for
$\Phi(h;x,y)$ as well as for $\Psi(h;x,y)$, but it was hitherto
unknown for nonlinear polynomials $h$ whether $\Phi(h;x,x^{g-\delta})$
even tended to infinity with $x$ for any fixed positive $\delta$. For
linear polynomials, Theorem \ref{primethm} is weaker in terms of the
admissible range of $\delta$ than existing theorems; for example,
Friedlander \cite{Fri:SPwLPF} shows that the lower bound
(\ref{primebd}) holds for any $\delta$ less than $1-1/(2\sqrt
e)=0.6967\dots$ when $h(t)=t+a$ for some nonzero integer $a$.

Finally, we establish by elementary means a theorem that in some
sense interpolates between Theorem \ref{mainthm} and Schinzel's result
(\ref{schin}):

\begin{theorem}
Let $h(t)$, $g$, and $k$ be as in Theorem \ref{mainthm}. Then when $x$
is sufficiently large, we have $\Psi(h;x,x^{g-1/k})\gg_h
x\log^{-k}x$. In particular, if $h$ is irreducible, then
$\Psi(h;x,x^{g-1})\gg_h x\log^{-1}x$.
\label{sillythm}
\end{theorem}
\noindent Theorem \ref{sillythm} has a weaker smoothness parameter
than (\ref{schin}) but provides a stronger quantitative lower bound,
while it has a stronger smoothness parameter than Theorem
\ref{mainthm} but a weaker lower bound.

Section \ref{grungysec} of this paper contains the outline of the
approach to establishing Theorems \ref{mainthm}--\ref{scantthm}, as
well as definitions of much of the notation used throughout the
paper. Section \ref{multsec} examines the multiplicative functions
that arise in the course of implementing this plan, while Section
\ref{mult2sec} addresses the asymptotics of sums of multiplicative
functions. Section \ref{sievesec} deals with the sieve-related work
and culminates in a proof of Proposition \ref{grungyprop}
below. Section \ref{primesec} provides an outline of the modifications
to this proof necessary to establish Theorem \ref{primethm}, and
Section \ref{sillysec} contains a proof of Theorem~\ref{sillythm}.

Throughout this paper, we use the usual notation $(m,n)$ and $[m,n]$
for the greatest common divisor and least common multiple,
respectively, of $m$ and $n$; $\mu(n)$ for the M\"obius function;
$\phi(n)$ for the Euler totient function; $d(n)$ for the number of
divisors of $n$; $\Lambda(n)$ for the von Mangoldt function; and
$\omega(n)$ for the number of distinct prime factors of $n$. We also
use the notation $m\mid n$ to mean that $m$ divides $n$, and
$p^r\exdiv n$ to mean that the prime power $p^r$ exactly divides $n$,
i.e., $p^r$ divides $n$ but $p^{r+1}$ does not. The constants implicit
in the $O$- and $\ll$ symbols in this paper may depend where
appropriate on the polynomial under investigation ($h$, when it
denotes a polynomial, or $f$) and on quantities defined only in terms
of that polynomial (e.g., $g$, $k$, $t_0$, and $\Delta$), and also on
$\delta$ and $\ep$; the same dependencies are allowed when the phrase
``sufficiently large'' is used.

The author would like to express his appreciation to Michael Bennett,
John Friedlander, Hugh Montgomery, Pieter Moree, Carl Pomerance, and
especially Trevor Wooley for helpful discussions regarding existing
results in this area and preliminary versions of this work. The author
would also like to thank Henryk Iwaniec for providing access to some
of his unpublished work and Jon Sorenson for providing translations of
several of the Russian papers cited herein. The author was supported
by a National Science Foundation Graduate Research Fellowship and by
National Science Foundation grant DMS 9304580.

%% file: grungy.tex
Theorems \ref{mainthm}--\ref{scantthm} will follow from the more
quantitative Proposition \ref{grungyprop} below. It is convenient to
define, for any integer-valued polynomial $h$, the quantities
$C(h;x,y)$ and $C(h;x,\sh,y)$, the complements of the quantities
$\Psi(h;x,y)$ and $\Psi(h;x,\sh,y)$:
\begin{equation*}
\begin{split}
C(h;x,y) &= x - \Psi(h;x,y), \\
C(h;x,\sh,y) &= \sh - \Psi(h;x,\sh,y).
\end{split}
\end{equation*}

\begin{proposition} Let $f(t)$ be an irreducible integer-valued
polynomial that is not identically zero modulo any prime. Let
$x\ge\sh\ge2$ and $0<\delta<1/2$ be real numbers, and set
$\xi=\max_{n\le x}|f(n)|$. Then
\begin{equation*}
C(f;x,\sh,\xi\sh^{-\delta}) \le \sh \big( {2\delta\over1-\delta} + O{(
\log^{-1/3}\sh )} \big).
\end{equation*}
\label{grungyprop}
\end{proposition}
\vskip-12pt

Let us see why Proposition \ref{grungyprop} implies Theorem
\ref{shortthm} for a general polynomial $h(t)$. First we let $m$ be
the largest integer such that $h$ is identically zero\mod m; i.e., $m$
is the greatest common divisor of the values $h(n)$ for $n\in\bf
Z$. If we set $h_1(t)=h(t)/m$, then $h_1$ is still integer-valued, and
furthermore $h_1$ is not identically zero modulo any prime by the
definition of $m$. Also, as long as $y$ is greater than $m$, the value
$h(n)$ is $y$-smooth precisely when $h_1(n)$ is $y$-smooth. Thus it
suffices to consider $\hht$.

Let $g$ be the largest degree of any irreducible factor of $\hht$
(equivalently, of $h$), and write
\begin{equation*}
\hht(t) = f_1(t)^{r_1}\dotsm f_k(t)^{r_k}\hhh(t),
\end{equation*}
where the $f_i$ are distinct irreducible polynomials of degree $g$
with integer coefficients and every irreducible factor of $\hhh$ has
degree at most $g-1$. Let $\xi_i=\max_{n\le x}|f_i(n)|$ and
$\xi=\max_i\xi_i$, and note that $\xi$ has order of magnitude $x^g$.
The values $\hhh(n)$ when $n\le x$ are always $O(x^{g-1})$-smooth; in
particular, given any $0<\delta<1$, they are $\xi\sh^{-\delta}$-smooth
for $x$ sufficiently large, since $L\le x$. Thus $\hht(n)$ fails
to be $\xi\sh^{-\delta}$-smooth precisely when at least one of the
$f_i(n)$ is not $\xi\sh^{-\delta}$-smooth, which implies that
\begin{equation}
C(\hht;x,\sh,\xi\sh^{-\delta}) \le \sum_{i=1}^k
C(f_i;x,\sh,\xi\sh^{-\delta}) \le \sum_{i=1}^k
C(f_i;x,\sh,\xi_i\sh^{-\delta}).
\label{separate}
\end{equation}
Since $\hht$ is not identically zero modulo any prime, the same is
true of each $f_i$. This allows us to apply Proposition
\ref{grungyprop} to each term in the latter sum in the inequality
(\ref{separate}), which becomes
\begin{equation}
C(\hht;x,\sh,\xi\sh^{-\delta}) \le k\sh \big( {2\delta \over1-\delta}
+ O(\log^{-1/3}\sh) \big).
\label{fishbd}
\end{equation}
Therefore $\Psi(h_1;x,\sh,\xi\sh^{-\delta}) = \sh -
C(h_1;x,\sh,\xi\sh^{-\delta}) \gg\sh$ whenever $\delta<(1+2k)^{-1}$,
which is the assertion of Theorem \ref{shortthm}, aside from the minor
difference between $\xi$ and $x^g$ which can be accommodated by a
very small change in $\delta$.

Theorem \ref{scantthm} follows from the inequality (\ref{fishbd}),
with $\sh=x$, whenever $\delta$ is small enough that $2\delta
k/(1-\delta)$ is less than the density $\eta$ of $\A$, which is
equivalent to the condition that $\delta<\eta/(2k+\eta)$. Theorem
\ref{mainthm} certainly follows from Proposition \ref{grungyprop} as
well, since it is the special case of Theorem \ref{shortthm} with
$\sh=x$, or a special case of Theorem \ref{scantthm} with $\eta=1$.

It is worth remarking that when $g=1$, Proposition \ref{grungyprop} is
a result about smooth integers in short intervals or short arithmetic
progressions. For the purposes of illustration, we take $h(t)$ to be
simply $t$, and put $L=x^\beta$ for some $0<\beta<1$ and
$\delta=(1-\alpha)/\beta$ for some $1-\beta/2<\alpha<1$, so that
$0<\delta<1/2$ and $\xi L^{-\delta}=x^\alpha$. Proposition
\ref{grungyprop} then gives us
\begin{equation}
\Psi(x,x^\alpha) - \Psi(x-x^\beta,x^\alpha) \ge \big(
{\beta-3(1-\alpha) \over \beta-(1-\alpha)} + o(1) \big) x^\beta,
\label{notFL}
\end{equation}
which is nontrivial in the range $\alpha+\beta/3>1$. Existing results
give nontrivial lower bounds for $\Psi(x,x^\alpha) -
\Psi(x-x^\beta,x^\alpha)$ for larger ranges of $\alpha$ and $\beta$
(see for instance Friedlander--Lagarias \cite[Theorem
2.4]{FriLag:OtDiSI}), so the lower bound (\ref{notFL}) is not
qualitatively new, although the constant
% $\big( \beta-3(1-\alpha) \big)/\big( \beta-(1-\alpha) \big)$
on the right-hand side seems to be an improvement over existing
results for certain values of $\alpha$ and $\beta$. Although
Proposition \ref{grungyprop} also gives an explicit lower bound for
the number of smooth integers in a short interval from a fixed
arithmetic progression, the methods in \cite{FriLag:OtDiSI} and
similar papers can surely be applied to this situation as well.

One can show that the values of $f$ are often free of small prime
factors by using a lower bound sieve to sieve out those values that
are multiples of small primes; however, this approach has no chance of
showing that the values of $f$ are often smooth if the degree of $f$
is at least 2. One would have to sieve by $\gg x^g$ primes,
necessitating a sum of $\gg x^g$ error terms. The most optimistic hope
would be that the individual error terms were uniformly bounded and
that we could obtain square-root cancellation in the sum of the error
terms, and even this would result in an error whose magnitude would be
$\gg x^{g/2}$, which would swamp the main term. Instead, we will
establish Proposition \ref{grungyprop} by bounding from above the
number of values of $f$ that are divisible by a prime greater than
$\xi\sh^{-\delta}$; broadly speaking, we will accomplish this by
grouping these values by their cofactors, the remainders when the
large prime divisors are removed from the values (equation
(\ref{Cf}) below contains an example of this grouping), and using
an upper bound sieve.

%% file: mult.tex
Let $f$ be an irreducible integer-valued polynomial that is not
identically zero modulo any prime, as in the statement of Proposition
\ref{grungyprop}. We define $\sigma(h)$ to be the number of solutions
of $f(x)\equiv0\mod h$. It is easily seen, by the Chinese remainder
theorem and the assumption that $f$ is not identically zero modulo any
prime, that $\sigma$ is a multiplicative function satisfying
$0\le\sigma(h)<h$. We also have a bound on $\sigma(h)$ in terms of the
degree $g$ and the discriminant $\Delta$ of $f$. Write $\Delta=\prod_p
p^{\nu(p)}$, where all but finitely many of the $\nu(p)$ are zero (the
discriminant $\Delta$ itself is nonzero because $f$ is irreducible).
Huxley \cite{Hux:ANoPC} gives a bound for $\sigma$ that implies
\begin{equation}
\sigma(p^r)\le gp^{\nu(p)/2}
\label{sandor}
\end{equation}
for any prime power $p^r$ (this estimate is improved by Stewart
\cite{Ste:OtNoSoPCaTE}, though it will suffice for our purposes as
stated). From the bound (\ref{sandor}) it follows that $\sigma(h)\le
g^{\omega(h)}\Delta^{1/2}\ll g^{\omega(h)}\ll h^\ep$ for any $\ep>0$,
since the implicit constants may depend on $f$ and $\ep$.

It is well-known that $\sigma(p)$ is equal to $1$ on average,
since $f$ is irreducible; in fact, Nagel \cite{Nag:GdTdT} showed that,
for any polynomial $H(t)$ with integer coefficients and with
$\sigma(H;p)$ roots\mod p for each prime $p$, the asymptotic formula
\begin{equation}
\sum_{p<w} \frac{\sigma(H;p)\log p}p = \kappa(H) \log w + O_H(1)
\label{nagel}
\end{equation}
holds for all $w\ge2$, where $\kappa(H)$ is the number of irreducible
factors of $H$. This readily implies that
\begin{equation}
\prod_{w_1\le p<w_2} \big( 1-\frac{\sigma(p)}p \big)^{-1} = \big(
{\log w_2\over\log w_1} \big) \big( 1 + O{\big( \frac1{\log w_1}
\big)} \big)
\label{rhoprodasym}
\end{equation}
for all $2\le w_1\le w_2$, or equivalently (by Mertens' formula)
\begin{equation}
\prod_{p<w} \big( 1-\frac{\sigma(p)}p \big)^{-1} = e^\gamma\log w
\big( 1 + O{\big( \frac1{\log w} \big)} \big) \prod_p \big(
1-\frac{\sigma(p)}p \big)^{-1} \big( 1-\frac1p \big)
\label{otherrhoprod}
\end{equation}
for all $w\ge2$, where $\gamma$ is Euler's constant.

We also define a multiplicative function $\sigma^*(h)$ by stipulating
that on prime powers $p^r$, we have
\begin{equation}
\sigma^*(p^r) = \sigma(p^r) - \frac{\sigma(p^{r+1})}p.
\label{sigstardef}
\end{equation}
We remark that $\sigma(p^{r+1})$ counts the number of roots of
$f\mod{p^{r+1}}$, each of which corresponds to a root of $f\mod{p^r}$
simply by reducing\mod{p^r}. Moreover, this correspondence is at most
$p$-to-1, i.e., $\sigma(p^{r+1}) \le p\sigma(p^r)$. Consequently,
$\sigma^*$ is a nonnegative function. Since $\sigma^*(p^r)$ obviously
does not exceed $\sigma(p^r)$ for any prime power $p^r$, we have that
$0\le\sigma^*(h)\le\sigma(h)$ for any $h$.

We also note that $\sigma^*(h)=0$ if and only if there is a prime $p$
dividing $h$ such that $\sigma(ph) = p\sigma(h)$, by the
multiplicativity of $\sigma$ and $\sigma^*$. This is equivalent to
saying that there is a prime $p$ such that every integer that is a
root of $f\mod h$ is also a root of $f\mod{ph}$; such a prime must
necessarily divide $h$, by the multiplicativity of $\sigma$ and the
assumption that $f$ is not identically zero modulo any prime.

Expressions of the form $\sigma(nh)/\sigma(h)$ will arise later in
connection with $\sigma^*$, and we will need to know that such
expressions are multiplicative in the variable $n$. This is a general
property of multiplicative functions which we establish in the
following lemma.

\begin{lemma}
If $g(n)$ is a multiplicative function, then for any fixed number $h$
satisfying $g(h)\ne0$, the function $g(nh)/g(h)$ is also a
multiplicative function of $n$.
\label{alsomult}
\end{lemma}

\begin{proof}
It is easily seen that a multiplicative function $g$ satisfies
$g(m)g(n) = g([m,n])g((m,n))$ for any numbers $m$ and $n$, by writing
all four arguments as products of prime powers. This implies that if
$n_1$ and $n_2$ are relatively prime, we have
\begin{equation*}
{g(n_1h)\over g(h)} {g(n_2h)\over g(h)} = {g([n_1h,n_2h])
g((n_1h,n_2h))\over g(h)^2} = {g(n_1n_2h)g(h) \over g(h)^2} =
{g(n_1n_2h)\over g(h)},
\end{equation*}
which establishes the lemma.
\end{proof}

We will also need the following upper bound when we apply the sieve in
Section \ref{sievesec}.

\begin{lemma}
For any positive integer $h$ such that $\sigma^*(h)>0$, and for any
real numbers $2\le w_1\le w_2$, we have
\begin{equation*}
\prod_{w_1<p\le w_2} \big( 1-{\sigma(ph)\over p\sigma(h)} \big)^{-1}
\le \big( {\log w_2\over\log w_1} \big) \big( 1+O\big( \frac1{\log
w_1} \big) \big),
\end{equation*}
where the implicit constant does not depend on $h$.
\label{secondlem}
\end{lemma}

\begin{proof}
We recall that $\sigma^*(h)\le\sigma(h)$, so that the assumption that
$\sigma^*(h)$ is positive implies that $\sigma(h)$ is also positive.
If a prime $p$ does not divide $h$, then $\sigma(ph) =
\sigma(p)\sigma(h)$ by the multiplicativity of $\sigma$. If $p$
divides $h$ but does not divide the discriminant $\Delta$ of $f$, then
every root $b$ of $f$\mod h must satisfy $f'(b)\not\equiv0$\mod p. In
this case, every root $b$ of $f$\mod h corresponds to exactly one root
of $f$\mod{ph} by Hensel's Lemma, and in particular,
$\sigma(ph)=\sigma(h)$ in this case. Therefore we can write
\begin{multline}
\prod_{w_1<p\le w_2} \big( 1-{\sigma(ph)\over p\sigma(h)} \big)^{-1} =
\prod_{w_1<p\le w_2} \big( 1-{\sigma(p)\over p} \big)^{-1} \prod
\begin{Sb}w_1<p\le w_2 \\ p\mid h,\, p\dnd\Delta\end{Sb} \big(
1-\frac{\sigma(p)}p \big) \big( 1-\frac1p \big)^{-1} \\
{}\times \prod \begin{Sb}w_1<p\le w_2 \\ p\mid h,\,
p\mid\Delta\end{Sb} \big( 1-\frac{\sigma(p)}p \big) \big(
1-{\sigma(ph)\over p\sigma(h)} \big)^{-1}.
\label{threeprinces}
\end{multline}

Equation (\ref{rhoprodasym}) gives an asymptotic formula for the first
product in this equation. Each term in the second product is at most
1, since the fact that $\sigma(h)>0$ certainly implies that $f$ has at
least one root\mod p for every prime $p$ dividing $h$, so that
$\sigma(p)\ge1$ for the primes in the second product of equation
(\ref{threeprinces}). Finally, the third product can be bounded
above by
\begin{equation*}
\prod \begin{Sb}w_1<p\le w_2 \\ p^r\exdiv h,\, p\mid\Delta\end{Sb}
\big( 1-{\sigma(p^{r+1}) \over p\sigma(p^r)} \big)^{-1}
\end{equation*}
by the fact that $\sigma$ is nonnegative and multiplicative.
Furthermore, since $\sigma^*(h)>0$ is equiv\-a\-lent to the condition
that $\sigma(p^{r+1}) < p\sigma(p^r)$ for every prime power $p^r$
exactly dividing $h$, this product can in turn be bounded by
\begin{equation*}
\prod \begin{Sb}p>w_1 \\ p\mid\Delta\end{Sb} \max \begin{Sb}r\ge0 \\
\sigma(p^{r+1}) < p\sigma(p^r)\end{Sb} \big( 1-{\sigma(p^{r+1}) \over
p\sigma(p^r)} \big)^{-1} \le \prod \begin{Sb}p>w_1 \\
p\mid\Delta\end{Sb} \max_{r\ge0} p\sigma(p^r) \le \prod
\begin{Sb}p>w_1 \\ p\mid\Delta\end{Sb} gp^{1+\nu(p)/2}
\end{equation*}
by the upper bound (\ref{sandor}).

Equation (\ref{threeprinces}) now becomes
\begin{equation*}
\prod_{w_1<p\le w_2} \big( 1-{\sigma(ph) \over p\sigma(h)} \big) \le
\big( {\log w_2\over\log w_1} \big) \big( 1+O\big( \frac1{\log w_1}
\big) \big) \prod \begin{Sb}p>w_1 \\ p\mid\Delta\end{Sb}
gp^{1+\nu(p)/2}.
\end{equation*}
This last product is bounded above independently of $h$, and it has
the value 1 as soon as $w_1$ exceeds $\Delta$. Therefore its
contribution can be absorbed into the implicit constant in the error
term. This establishes the lemma.
\end{proof}

%% file: mult2.tex
Our primary goal for this section is to establish an asymptotic
formula for a summatory function $M_g(x)$ associated with a
multiplicative function $g(n)$, defined by
\begin{equation*}
M_g(x) = \sum_{n\le x} \gov n.  \label{Mgdef}
\end{equation*}
We are interested in an asymptotic formula for $M_g(x)$ when $g(p)$ is
constant on average over primes, as is usually the case for the
multiplicative functions that arise in sieve problems. Specifically,
we impose the condition that there is a constant $\kappa=\kappa(g)$
such that
\begin{equation}
\sum_{p\le x} \glov p = \kappa\log x + O_g(1)  \label{consteq}
\end{equation}
for all $x\ge2$.

Although the ideas used in establishing the following proposition have
been part of the ``folklore'' for some time, the literature does not
seem to contain a result in precisely this form. Wirsing's pioneering
work \cite{Wir:DAVvS}, for instance, requires $g$ to be a nonnegative
function and implies an asymptotic formula for $M_g(x)$ without a
quantitative error term; while Halberstam and Richert \cite[Lemma
5.4]{HalRic:SM} give an analogous result with a quantitative error
term, but one that requires $g$ to be supported on squarefree integers
in addition to being nonnegative. Both results are slightly too
restrictive for our purposes as stated.

Consequently we provide a self-contained proof of an asymptotic
formula for $M_g(x)$ with a quantitative error term, for
multiplicative functions $g$ that are not necessarily supported on
squarefree integers. The proof below, which is based on unpublished
work of Iwaniec (used with his kind permission) that stems from ideas
of Wirsing and Chebyshev, has the advantage that $g$ is freed from the
requirement of being nonnegative. We state the result in a more
general form than is required for our present purposes, with a mind
towards other applications and because the proof is exactly the same
in the more general setting.

\begin{proposition}
Suppose that $g(n)$ is a complex-valued multiplicative function such
that the asymptotic formula {\rm(\ref{consteq})} holds for some
complex number $\kappa=\xi+i\eta$ satisfying $\eta^2<2\xi+1$ (so that
$\xi>-1/2$ in particular). Suppose also that
\begin{equation}
\sum_p \aglov p \sum_{r=1}^\infty \agov{p^r} + \sum_p
\sum_{r=2}^\infty \aglov{p^r} < \infty,
\label{minorconverge}
\end{equation}
and that there exists a nonnegative real number $\beta=\beta(g)<\xi+1$
such that
\begin{equation}
\prod_{p\le x} \big( 1 + \agov{p} \big) \ll_g \log^{\beta} x
\label{minorprod}
\end{equation}
for all $x\ge2$. Then the asymptotic formula
\begin{equation}
M_g(x) = c(g)\log^\kappa x + O_g((\log x)^{\beta-1})  \label{mvform}
\end{equation}
holds for all $x\ge2$, where $\log^\kappa x$ denotes the principal
branch of $t^\kappa$, and $c(g)$ is defined by the convergent product
\begin{equation}
c(g) = \Gamma(\kappa+1)^{-1} \prod_p \big( 1-\frac1p \big)^\kappa \big(
1 + \gov p + \gov{p^2} + \dotsb \big).
\label{cgdef}
\end{equation}
\label{mvprop}
\end{proposition}

\noindent We remark that the condition (\ref{minorprod}) cannot hold
with any $\beta<\abs\kappa$ if $g$ satisfies the asymptotic formula
(\ref{consteq}). The necessity that $\beta$ be less than $\xi+1$, so
that the formula (\ref{mvform}) is truly an asymptotic formula,
requires us to consider only those $\kappa$ for which
$\abs\kappa<\xi+1$; this is the source of the condition
$\eta^2<2\xi+1$ on $\kappa$.

The conditions (\ref{minorconverge}) and (\ref{minorprod}) are usually
very easily verified in practice. For example, the condition
(\ref{minorconverge}) automatically holds if there is a constant
$\alpha<1/2$ such that $g(n)\ll n^\alpha$; and if $g$ is in fact a
nonnegative function (so that, in particular, $\kappa$ is
nonnegative), then the condition (\ref{minorprod}), with
$\beta=\kappa$, follows from the asymptotic formula
(\ref{consteq}). We also remark that from equation (\ref{mvform}), it
follows easily by partial summation that
\begin{equation}
\sum_{n<x}g(n)\ll_g x\log^{\beta-1}x  \label{mvcor}
\end{equation}
under the hypotheses of the proposition.

\vskip12pt

\begin{proof}
All of the constants implicit in the $O$- and $\ll$ symbols in this
proof may depend on the multiplicative function $g$, and thus on
$\kappa$ and $\beta$ as well. We begin by examining an analogue of
$M_g(x)$ weighted by a logarithmic factor. We have
\begin{equation}
\begin{split}
\sum_{n\le x} \glov n &= \sum_{n\le x} \gov n \sum_{p^r\mid\mid n}
\log p^r \\
&= \sum_{r=1}^\infty \sum_{p\le x^{1/r}} \glov{p^r} \sum
\begin{Sb}m\le x/p^r \\ p\dnd m\end{Sb} \gov m \\
&= \sum_{p\le x} \glov p \sum_{m\le x/p} \gov m - \sum_{p\le x} \glov
p \sum \begin{Sb}m\le x/p \\ p\mid m\end{Sb} \gov m \\
& \qquad {}+ \sum_{r=2}^\infty \sum_{p\le x^{1/r}} \glov{p^r} \sum
\begin{Sb}m\le x/p^r \\ p\dnd m\end{Sb} \gov m \\
&= \Sigma_1 - \Sigma_2 + \Sigma_3,
\end{split}
\label{sigmaseq}
\end{equation}
say. If we define the function $\delta(x)$ by
\begin{equation}
\delta(x) = \sum_{p\le x} \glov p - \kappa\log x,  \label{deltadef}
\end{equation}
then $\Sigma_1$ becomes
\begin{equation}
\Sigma_1 = \sum_{m\le x} \gov m \sum_{p\le x/m} \glov p = \kappa
\sum_{m\le x} \gov m \log\frac xm + \sum_{m\le x} \gov m \delta\big(
\frac xm \big).
\label{using}
\end{equation}
Since $M_g(x)=1$ for $1\le x<2$ and
\begin{equation*}
M_g(x)\log x - \sum_{m\le x} \frac{g(m)\log m}m = \sum_{m\le x} \gov
m\log\frac xm = \int_1^x M_g(t) \frac{dt}t
\end{equation*}
by partial summation, we can rewrite equation (\ref{sigmaseq}) using
equation (\ref{using}) as
\begin{equation}
M_g(x)\log x - (\kappa+1) \int_2^x M_g(t)t^{-1} dt = E_g(x),
\label{usetwice}
\end{equation}
where we have defined
\begin{equation}
E_g(x) = (\kappa+1)\log 2 + \sum_{m\le x} \gov m \delta\big( \frac xm \big)
- \Sigma_2 + \Sigma_3.  \label{Egdef}
\end{equation}

We integrate both sides of equation (\ref{usetwice}) against
$x^{-1}(\log x)^{-\kappa-2}$, obtaining
\begin{multline}
\int_2^x M_g(u)u^{-1}(\log u)^{-\kappa-1} du - (\kappa+1)\int_2^x u^{-1}(\log
u)^{-\kappa-2} \int_2^u M_g(t) t^{-1} dt\,du \\
= \int_2^x E_g(u) u^{-1}(\log u)^{-\kappa-2} du.
\label{okthen}
\end{multline}
Some cancellation can be obtained on the left-hand side by switching
the order of integration in the double integral and evaluating the new
inner integral; equation (\ref{okthen}) becomes simply
\begin{equation*}
(\log x)^{-\kappa-1} \int_2^x M_g(u)u^{-1} du = \int_2^x E_g(u) u^{-1}(\log
u)^{-\kappa-2} du.
\end{equation*}
We can substitute this into equation (\ref{usetwice}), divide by $\log
x$, and rearrange terms to get
\begin{equation}
M_g(x) = (\kappa+1) \log^\kappa x \int_2^x E_g(u) u^{-1}(\log
u)^{-\kappa-2} du + E_g(x)\log^{-1}x.
\label{rearrange}
\end{equation}

An upper bound for $E_g(x)$ is now needed. Since $\delta(x)$ is
bounded from its definition (\ref{deltadef}) and the asymptotic
formula (\ref{consteq}), we have
\begin{equation}
\sum_{m\le x} \gov m \delta\big( \frac xm \big) \ll \sum_{m\le x}
\agov m.
\label{modpf}
\end{equation}
We also have
\begin{equation}
\sum_{m\le x} \agov m \le \prod_{p\le x} \bigg( 1 + \sum_{r=1}^\infty
\agov{p^r} \bigg) \le \prod_{p\le x} \big( 1+\agov p \big) \prod_{p\le
x} \bigg( 1 + \sum_{r=2}^\infty \agov{p^r} \bigg).
\label{protoMabsg}
\end{equation}
Because the sum $\sum_p \sum_{r=2}^\infty \abs{g(p^r)}/p^r$ converges
by the hypothesis (\ref{minorconverge}), the last product in equation
(\ref{protoMabsg}) is bounded as $x$ tends to infinity. Therefore the
hypothesis (\ref{minorprod}) implies that
\begin{equation}
\sum_{m\le x} \agov m \ll \log^{\beta}x.
\label{Mabsg}
\end{equation}
The terms $\Sigma_2$ and $\Sigma_3$ can be estimated by
\begin{equation*}
\Sigma_2 = \sum_{p\le x} \glov p \sum_{r=1}^\infty \gov{p^r} \sum
\begin{Sb}l\le x/p^{r+1} \\ p\dnd l\end{Sb} \gov l \ll \sum_{p\le x}
\aglov p \sum_{r=1}^\infty \agov{p^r} \sum_{l\le x} \agov l
\end{equation*}
and
\begin{equation*}
\Sigma_3 \ll \sum_{p\le x} \sum_{r=2}^\infty \aglov{p^r} \sum_{m\le
x} \agov m,
\end{equation*}
and so both $\Sigma_2$ and $\Sigma_3$ are $\ll \log^{\beta}x$ by the
estimate (\ref{Mabsg}) and the hypothesis (\ref{minorconverge}).
Therefore, by the definition (\ref{Egdef}) of $E_g(x)$, we see that
\begin{equation}
E_g(x) \ll \log^{\beta}x.  \label{Egxbd}
\end{equation}

In particular, since $\beta<\xi+1$, we have
\begin{equation}
\int_x^\infty E_g(u) u^{-1}(\log u)^{-\kappa-2} du \ll \int_x^\infty
u^{-1}(\log u)^{\beta-\xi-2} du \ll (\log x)^{\beta-\xi-1},
\label{tenw}
\end{equation}
and so equation (\ref{rearrange}) and the bound (\ref{Egxbd}) give us
the asymptotic formula
\begin{equation}
M_g(x) = c(g) \log^\kappa x + O((\log x)^{\beta-1})  \label{allbutcg}
\end{equation}
for $x\ge2$, where
\begin{equation}
c(g) = (\kappa+1) \int_2^\infty E_g(u) u^{-1}(\log u)^{-\kappa-2} du.
\end{equation}

To complete the proof of the proposition, we need to show that $c(g)$
can be written in the form given by (\ref{cgdef}); we accomplish this
indirectly, using the asymptotic formula (\ref{allbutcg}). Consider
the zeta-function $\zeta_g(s)$ formed from $g$, defined by
\begin{equation*}
\zeta_g(s) = \sum_{n=1}^\infty {g(n)\over n^s}.
\end{equation*}
From the estimate (\ref{Mabsg}) and partial summation, we see that
$\zeta_g(s)$ converges absolutely for $s>1$ (we will only need to
consider real values of $s$), and thus has an Euler product
representation
\begin{equation}
\zeta_g(s) = \prod_p \big( 1 + {g(p)\over p^s} + {g(p^2)\over p^{2s}}
+ \dotsb \big)
\label{eulerg}
\end{equation}
for $s>1$.

We can also use partial summation to write
\begin{equation}
\zeta_g(s+1) = s \int_1^\infty M_g(t) t^{-s-1}dt
\label{homestretch}
\end{equation}
for $s>0$. Since $M_g(x)=1$ for $1\le x<2$, it is certainly true that
\begin{equation*}
M_g(x) = c(g) \log^\kappa x + O(1+\log^\xi x)
\end{equation*}
in that range; using this together with the asymptotic formula
(\ref{allbutcg}), equation (\ref{homestretch}) becomes
\begin{multline*}
\zeta_g(s+1) = s \int_1^\infty c(g) \log^\kappa t \cdot t^{-s-1}
dt \\
{}+ O\bigg( s \int_1^2 (1+\log^\xi t) t^{-s-1}dt + s \int_2^\infty
(\log t)^{\beta-1} t^{-s-1}dt \bigg),
\end{multline*}
valid uniformly for $s>0$. Making the change of variables $t=e^{u/s}$
in all three integrals and multiplying through by $s^\kappa$ yields
\begin{equation}
\begin{split}
s^\kappa \zeta_g(s+1) &= c(g) \int_0^\infty u^\kappa e^{-u}du +
O\bigg( \int_0^{s\log2} (s^\xi+u^\xi)e^{-u}du + s^{\xi-\beta+1}
\int_{s\log2}^\infty u^{\beta-1}e^{-u}du \bigg) \\
&= c(g)\Gamma(\kappa+1) + O(s^{\xi-\beta+1}\log s^{-1})
\end{split}
\label{sto}
\end{equation}
as $s\to0^+$, where the exponent $\xi-\beta+1$ is positive and at most
1 (since $\beta\ge\abs\kappa\ge\xi$). Because the Riemann
$\zeta$-function satisfies $s\zeta(s+1)=1+O(s)$ as $s\to0^+$, equation
(\ref{sto}) implies
\begin{equation}
\zeta(s+1)^{-\kappa} \zeta_g(s+1) = c(g)\Gamma(\kappa+1) +
O(s^{\xi-\beta+1}\log s^{-1}).
\label{recent}
\end{equation}

On the other hand, from equation (\ref{eulerg}) we certainly have the
Euler product representation
\begin{equation*}
\zeta(s+1)^{-\kappa} \zeta_g(s+1) = \prod_p \big( 1-\frac1{p^{s+1}}
\big)^\kappa \big( 1 + {g(p)\over p^{s+1}} + {g(p^2)\over p^{2(s+1)}}
+ \dotsb \big)
\end{equation*}
for $s>0$, and one can show that in fact this Euler product converges
uniformly for $s\ge0$. The important contribution comes from the sum
$\sum_p (g(p)-\kappa)/p^{s+1}$, and we see from the hypothesis
(\ref{consteq}) and partial summation that
\begin{equation*}
\sum_{p>x} {g(p)-\kappa\over p^{s+1}} \ll \frac1{x^s\log x}
\end{equation*}
uniformly for $s\ge0$ and $x\ge2$. The remaining contributions can be
controlled using the hypothesis (\ref{minorconverge}).

Consequently, taking the limit of both sides of equation
(\ref{recent}) as $s\to0^+$ gives us
\begin{equation*}
\prod_p \big( 1-\frac1{p} \big)^\kappa \big( 1 + {g(p)\over p} +
{g(p^2)\over p^{2}} + \dotsb \big) = c(g) \Gamma(\kappa+1)
\end{equation*}
(where we have just shown that the product on the left-hand side
converges), which is equivalent to (\ref{cgdef}). This establishes the
proposition.
\end{proof}

The following proposition gives a similar asymptotic formula for the
restricted sum
\begin{equation*}
M_g(x,q) = \sum \begin{Sb}n\le x \\ (n,q)=1\end{Sb} \gov n.
\end{equation*}
Although we will not need such a formula in this paper, results of
this type have widespread applicability, and so we include it also
with a mind towards other applications.

\begin{proposition}
Suppose that $g(n)$ satisfies the hypotheses of Proposition
\ref{mvprop}. Then the asymptotic formula
\begin{equation*}
M_g(x,q) = c_q(g)\log^\kappa x + O_g(\delta(q)(\log x)^{\beta-1})
\end{equation*}
holds uniformly for all $x\ge2$ and all nonzero integers $q$, where
\begin{equation*}
c_q(g) = \Gamma(\kappa+1)^{-1} \big( \frac{\phi(q)}q \big)^\kappa
\prod_{p\dnd q} \big( 1-\frac1p \big)^\kappa \big( 1 + \gov p +
\gov{p^2} + \dotsb \big)
\end{equation*}
and $\delta(q) = 1 + \sum_{p\mid q} \abs{g(p)}(\log p)/p$.
\label{coprimeprop}
\end{proposition}

\noindent We remark that we can also write
\begin{equation*}
c_q(g) = c(q) \prod_{p\mid q} \big( 1 + \gov p + \gov{p^2} + \dotsb
\big)^{-1},
\end{equation*}
as long as no term $(1+g(p)/p+g(p^2)/p^2+\dotsb)$ sums to zero.

\vskip12pt

\begin{proof}
We would like to apply Proposition \ref{mvprop} to the multiplicative
function $g_q(n)$ defined by
\begin{equation*}
g_q(n) = \begin{cases}g(n), &\rmif (n,q)=1, \\ 0, &\rmif
(n,q)>1.\end{cases}
\end{equation*}
Certainly $\abs{g_q(n)}\le\abs{g(n)}$, and so the estimates
(\ref{minorconverge}) and (\ref{minorprod}) for $g_q$ follow from the
same estimates for $g$. We also have
\begin{equation*}
\begin{split}
\sum_{p\le x} \frac{g_q(p)\log p}p = \sum \begin{Sb}p\le x \\ p\dnd
q\end{Sb} \glov p &= \sum_{p\le x} \glov p - \sum \begin{Sb}p\le x \\
p\mid q\end{Sb} \glov p \\
&= \kappa\log x + O_g(1) + O\bigg( \sum_{p\mid q} \aglov p \bigg)
\end{split}
\end{equation*}
from the assumption that $g$ satisfies equation (\ref{consteq}).
Therefore $g_q$ satisfies equation (\ref{consteq}) as well, with the
error term being $\ll_g \delta(q)$ uniformly in $x$.

If we keep this dependence on $q$ explicit throughout the proof of
Proposition \ref{mvprop}, the only modification necessary is to
include a factor of $\delta(q)$ on the right-hand sides of the
estimates (\ref{modpf}), (\ref{Egxbd}), and (\ref{tenw}) and in the
error term in equation (\ref{allbutcg}). Therefore, the application of
Proposition \ref{mvprop} to $g_q$ yields
\begin{equation*}
M_g(x,q) = M_{g_q}(x) = c(g_q) \log^\kappa x + O_g(\delta(q)(\log
x)^{\beta-1}),
\end{equation*}
where the implicit constant is independent of $q$. Because
\begin{equation*}
\begin{split}
c(g_q) &= \Gamma(\kappa+1)^{-1} \prod_p \big( 1-\frac1p \big)^\kappa \big(
1 + \gqov p + \gqov{p^2} + \dotsb \big) \\
&= \Gamma(\kappa+1)^{-1} \prod_{p\mid q} \big( 1-\frac1p \big)^\kappa
\prod_{p\dnd q} \big( 1-\frac1p \big)^\kappa \big( 1 + \gov p +
\gov{p^2} + \dotsb \big) = c_q(g),
\end{split}
\end{equation*}
the proposition is established.
\end{proof}

%% file: Sd.tex
In this section we reformulate Proposition \ref{grungyprop} in a way
that makes it amenable to treatment by sieve techniques. As in the
statement of Proposition \ref{grungyprop}, we let $x\ge\sh\ge2$ and
$0<\delta<1/2$ be real numbers. We can multiply $f$ by $-1$ if
necessary to make the leading coefficient positive without affecting
the smoothness of the values; and we can replace $f(t)$ by $f(t+t_0)$
for any fixed $t_0$ depending on $f$, since this only changes
$C(f;x,L,\xi L^{-\delta})$ by $O(1)$. Thus we may assume without loss
of generality that $f(n)$ is positive when $n$ is positive.

Put $\xi=\max_{x-L<n\le x}f(n)$, and notice that $L^\delta < \xi^{1/2}
< \xi L^{-\delta}$ when $x$ is sufficiently large, since $\xi\gg
x^g\ge x$ and $\delta<\half$. Letting $p$ denote only primes, we have
\begin{equation}
\begin{split}
C(f;x,L,\xi L^{-\delta}) &= \3{x-L<n\le x: \exists p>\xi L^{-\delta}
\hbox{ such that } p\mid f(n)} \\
&= \3{(n,p,h): x-L<n\le x,\, p>\xi L^{-\delta},\, f(n)=ph} \\
&= \sum_{h\ge1} \3{(n,p): x-L<n\le x,\, p>\xi L^{-\delta},\,
f(n)=ph}.
\end{split}
\label{Cf}
\end{equation}
(The integer $h$ plays the role of the cofactor mentioned at the end
of Section \ref{grungysec}.)

It is clear that $h$ must not exceed $L^\delta$ if it is to contribute
to this sum. Moreover, we claim that only those $h$ for which
$\sigma^*(h)>0$ contribute to the sum. If $\sigma^*(h)=0$, then by the
remarks following the definition (\ref{sigstardef}) of $\sigma^*$,
there is a prime $q$ dividing $h$ such that whenever $f(n)$ is
divisible by $h$, it is also the case that $f(n)/h$ is divisible by
$q$. But $q\le h<L^\delta < \xi L^{-\delta}$, and so there are no
pairs $(n,p)$ satisfying the description on the last line of equation
(\ref{Cf}). We may therefore write
\begin{equation}
C(f;x,L,\xi L^{-\delta}) = \sum \begin{Sb}1\le h<L^\delta \\
\sigma^*(h)>0\end{Sb} \3{(n,p): x-L<n\le x,\, p>\xi L^{-\delta},\,
f(n)=ph}.
\label{initial2}
\end{equation}
We remark that the purpose of insisting on this addition condition is
to facilitate the passage from equation (\ref{unstard}) to equation
(\ref{usedsieve}) in the proof of Lemma \ref{cancellem} below. If we
retained those terms for which $\sigma^*(h)=0$, a formal use of an
upper bound sieve and a mean value theorem for multiplicative
functions would result in infinite products containing local factors
equaling zero and infinity (respectively). However, these factors
would formally cancel at the end of the proof of Lemma \ref{cancellem}
along with the rest of the local factors, and so we see that the
restriction is technical rather than substantive.

To estimate the right-hand side of equation (\ref{initial2}) using an
upper bound sieve, we replace occurrences of the prime $p$ by any
integer $m$ whose prime factors are large. We define
\begin{equation}
S(z) = \sum \begin{Sb}1\le h<L^\delta \\ \sigma^*(h)>0\end{Sb}
\3{(n,m): x-L<n\le x,\, f(n)=mh,\, p\mid m\implies p>z}
\label{Szdef}
\end{equation}
and notice that the right-hand side of equation (\ref{initial2}) is
precisely $S(\xi\sh^{-\delta})$. It is clear that $S(z)$ is a
decreasing function of $z$, and therefore {\it to establish
Proposition \ref{grungyprop} and thus Theorems
\ref{mainthm}--\ref{scantthm}, it suffices to show that
\begin{equation}
S(z) \le \sh \big( {2\delta\over1-\delta} + O{( \log^{-1/3}\sh )} \big)
\label{italics}
\end{equation}
for some value of $z$ in the range $2\le z\le \xi\sh^{-\delta}$.\/}

As is standard in sieve problems, to understand $S(z)$ we need to
understand the corresponding sums over multiples of a given integer,
and so we define
\begin{equation}
S_d = \sum \begin{Sb}1\le h<L^\delta \\ \sigma^*(h)>0\end{Sb}
\3{(n,m): x-L<n\le x,\, f(n)=mh,\, d\mid m}.
\label{Sddef}
\end{equation}
We see that
\begin{equation}
\begin{split}
S_d &= \sum \begin{Sb}1\le h<L^\delta \\ \sigma^*(h)>0\end{Sb}
\3{x-L<n\le x: dh\mid f(n)} \\
&= \sum \begin{Sb}1\le h<L^\delta \\ \sigma^*(h)>0\end{Sb} \big(
{L\sigma(dh)\over dh} + O(\sigma(dh)) \big) = \frac Ld \sum
\begin{Sb}1\le h<L^\delta \\ \sigma^*(h)>0\end{Sb} \frac{\sigma(dh)}h
+ O\bigg( \sum_{1\le h<L^\delta} \sigma(dh) \bigg),
\end{split}
\label{sdeasy}
\end{equation}
since every block of $dh$ consecutive integers contains precisely
$\sigma(dh)$ roots of $f\mod{dh}$. We remark that we could evaluate
the sums over $h$ asymptotically at this time by
Proposition~\ref{mvprop}, but as those familiar with sieve methods
will recognize, it is crucial to keep the error term in the formula
for $S_d$ as small as possible before the sieve is applied. The use of
Proposition~\ref{mvprop} would permit only a relative error of
$\log^{-1}\sh$ in this formula, which would never allow us to sieve by
a set of primes up to a power of~$\sh$.

%% file: sieve.tex
We now describe the upper bound linear sieve introduced by Rosser and
developed by Iwaniec \cite{Iwa:ANFotETitLS, Iwa:RS}. For a real number
$w\ge2$, let $P(w)=\prod_{p<w}p$. Let $D>1$ be a real number, and
define a sequence $\{\lmp d\}$ of real numbers that are supported on
the squarefree numbers not exceeding $D$ as follows: let $\lmp
1=1$, and if $d=p_1\cdots p_r$ with $p_1>\cdots>p_r$, define
\begin{equation*}
\lmp d = \begin{cases}
(-1)^r, & \hbox{if $p_1\cdots p_{2i}p_{2i+1}^3<D$ for all $0\le
i<r/2$,} \\
0, & \hbox{otherwise.}
\end{cases}
\end{equation*}

\def\mm{M}
The sequence $\{\lmp d\}$ is in fact an upper bound sieve, that is, it
satisfies
\begin{equation}
\sum_{d\mid n} \lmp d \ge \sum_{d\mid n} \mu(d) = \begin{cases} 1,
&\rmif n=1, \\ 0, &\rmif n>1, \end{cases}
\label{upper}
\end{equation}
where the latter equality is the characteristic property of the
M\"obius function. In addition, Iwaniec \cite[Lemma
3]{Iwa:ANFotETitLS} shows that, uniformly for all multiplicative
functions $\mm$ satisfying $0\le \mm(p)<p$ for all primes $p$ and
\begin{equation}
\prod_{w_1\le p<w_2} \big( 1-\frac{\mm(p)}p \big)^{-1} \le \big( {\log
w_2\over\log w_1} \big) \big( 1 + O\big( \frac 1{\log w_1} \big) \big)
\label{twoconds}
\end{equation}
for all $2\le w_1\le w_2$, we have
\begin{equation}
\sum_{d\mid P(z)} \lmp d\frac{\mm(d)}d \le \prod_{p<z} \big(
1-\frac{\mm(p)}p \big) (F(s) + O(\log^{-1/3}D) ),
\label{upperlinear}
\end{equation}
for all $2\le z\le D$, where $s=\log D/\log z$. Here $F(u)$ is the
traditional upper-bound function of the linear sieve: it is the
continuous solution for $u>0$ of the system of differential-difference
equations \def\ff{f}
\begin{equation}
\begin{split}
F(u)=\frac{2e^\gamma}u \quad\hbox{and}\quad \ff(u)=0 \qquad
&(0<u\le2), \\
(uF(u))' = \ff(u-1) \quad\hbox{and}\quad (u\ff(u))' = F(u-1) \qquad
&(2<u).
\end{split}
\label{Ffdef}
\end{equation}
One can see that $F(u)$ and $\ff(u)$ are both nonnegative functions
and hence that $uF(u)$ is nondecreasing. (Of course, the companion
function $\ff(u)$ to $F(u)$ is not the same as the polynomial $f$
whose values we are investigating; we will not need to refer to this
companion function again, so no confusion should arise.)

%% file: cancel.tex
With this notation in place, we can provide an upper bound for an
expression that will arise in the main term of our sieve estimate for
$S(z)$.

\begin{lemma}
Let $\sh\ge2$ and $0<\delta<1/2$ be real numbers. For any real numbers
$z$ and $D$ satisfying $L^\delta\le z\le D\le\exp(\log^3L)$, we have
\begin{equation}
\sum_{d\mid P(z)} \frac{\lmp d}d \sum \begin{Sb}1\le h<L^\delta \\
\sigma^*(h)>0\end{Sb} {\sigma(dh)\over h} \le \big(
{2\delta\log\sh\over\log D} \big) \big( {sF(s)\over2e^\gamma} +
O(s\log^{-1/3}D) \big),
\label{canceleq}
\end{equation}
where $s=\log D/\log z$.
\label{cancellem}
\end{lemma}
\noindent We remark that the right-hand side of the inequality
(\ref{canceleq}) has no local factors depending on the polynomial
$f$. This should not be surprising, as the upper bound sieve
$\lambda_d$ is meant to mimic the behavior of $\mu(d)$, so that the
sum on the left-hand side of (\ref{canceleq}) should behave like
$\sum_d \sum_h \mu(d) \sigma(dh)/dh$. But for any multiplicative
function $M$, we formally have
\begin{equation*}
\sum_d \sum_h \mu(d) M(dh) = \sum_n M(n) \sum_{d\mid n} \mu(d) = M(1)
= 1.
\end{equation*}
We also remark that we have collected the terms in the upper bound
(\ref{canceleq}) in such a way as to highlight the quantity
$sF(s)/2e^\gamma$. Since $sF(s)$ is a nondecreasing function, as noted
after equation (\ref{Ffdef}), we should take $s$ to be as small as
possible (subject to $z\le D$) when applying the asymptotic inequality
(\ref{canceleq}).  Thus we set $D=z$, whence $s=1$ and
$sF(s)/2e^\gamma=1$ as well, again by (\ref{Ffdef}). We obtain
\begin{equation}
\sum_{d\mid P(z)} \frac{\lmp d}d \sum \begin{Sb}1\le h<L^\delta \\
\sigma^*(h)>0\end{Sb} {\sigma(dh)\over h} \le \big(
{2\delta\log\sh\over\log z} \big) (1+O(\log^{-1/3}z))
\label{Disz}
\end{equation}
for any $L\ge2$ and $z\ge L^\delta$.
\vskip12pt

\begin{proof}
We begin by recalling that $\sigma^*(h)\le\sigma(h)$, and so
$\sigma(h)$ is positive whenever $\sigma^*(h)$ is positive. With this
observation, we may write
\begin{equation}
\sum_{d\mid P(z)} \frac{\lmp d}d \sum \begin{Sb}1\le h<L^\delta \\
\sigma^*(h)>0\end{Sb} {\sigma(dh)\over h} = \sum \begin{Sb}1\le
h<L^\delta \\ \sigma^*(h)>0\end{Sb} \frac{\sigma(h)}h \sum_{d\mid
P(z)} \lmp d {\sigma(dh)\over d\sigma(h)}.
\label{unstard}
\end{equation}

By Lemma \ref{alsomult}, the function $M(d) = \sigma(dh)/\sigma(h)$ is
a multiplicative function of $d$, and so we would like to apply the
upper bound (\ref{upperlinear}) to the inner sums on the right-hand
side of equation (\ref{unstard}). The inequality (\ref{twoconds}) is
satisfied uniformly in $h$ by Lemma \ref{secondlem}, and so it
remains only to verify that $M(p)<p$ for every prime $p$. But if this
were not the case, then we would have a prime $p$ for which
$\sigma(ph)=p\sigma(h)$. By the comments following the definition
(\ref{sigstardef}) of $\sigma^*$, this would then imply that
$\sigma^*(h)=0$, and these values of $h$ are excluded from the sum in
equation (\ref{unstard}).

We are therefore allowed to apply the upper bound (\ref{upperlinear}),
and the resulting error term will be uniform in $h$ as well. Equation
(\ref{unstard}) thus becomes
\begin{equation}
\begin{split}
\sum_{d\mid P(z)} \frac{\lmp d}d \sum \begin{Sb}1\le h<L^\delta \\
\sigma^*(h)>0\end{Sb} {\sigma(dh)\over h} &\le \sum \begin{Sb}1\le
h<L^\delta \\ \sigma^*(h)>0\end{Sb} \frac{\sigma(h)}h \prod_{p<z}
\big( 1-\frac{\sigma(ph)}{p\sigma(h)} \big) (F(s)+O(\log^{-1/3}D)) \\
&= (F(s)+O(\log^{-1/3}D)) \prod_{p<z} \big( 1-\frac{\sigma(p)}p \big)
\\
&\qquad{}\times \sum \begin{Sb}1\le h<L^\delta \\ \sigma^*(h)>0\end{Sb}
\frac{\sigma(h)}h \prod \begin{Sb}p<z \\ p^r\exdiv h\end{Sb} \big(
1-\frac{\sigma(p)}p \big)^{-1} \big( 1-{\sigma(p^{r+1}) \over
p\sigma(p^r)} \big)
\end{split}
\label{usedsieve}
\end{equation}
by the multiplicativity of $\sigma$. Equation (\ref{otherrhoprod})
immediately gives the asymptotic formula
\begin{equation}
\prod_{p<z} \big( 1-\frac{\sigma(p)}p \big) = \frac1{e^\gamma\log z}
\big( 1 + O\big( \frac1{\log z} \big) \big) \prod_p \big(
1-\frac{\sigma(p)}p \big) \big( 1-\frac1p \big)^{-1}
\label{bunga}
\end{equation}
for the first product in the last expression of inequality
(\ref{usedsieve}).

If we define a multiplicative function $G(h)$ by $G(h) = \prod_{p\mid
h} (1-\sigma(p)/p)^{-1}$ and use the assumption that $z\ge L^\delta$,
then the sum over $h$ in the last line of inequality (\ref{usedsieve})
becomes
\begin{equation}
\begin{split}
\sum \begin{Sb}1\le h<L^\delta \\ \sigma^*(h)>0\end{Sb}
\frac{\sigma(h)}h \prod \begin{Sb}p<z \\ p^r\exdiv h\end{Sb} \big(
1-\frac{\sigma(p)}p \big)^{-1} \big( 1-{\sigma(p^{r+1}) \over
p\sigma(p^r)} \big) &= \sum \begin{Sb}1\le h<L^\delta \\
\sigma^*(h)>0\end{Sb} \frac1h G(h) \prod_{p^r\exdiv h} \sigma(p^r)
\big( 1-{\sigma(p^{r+1}) \over p\sigma(p^r)} \big) \\
&= \sum \begin{Sb}1\le h<L^\delta \\ \sigma^*(h)>0\end{Sb}
\frac{G(h)\sigma^*(h)}h
\end{split}
\label{goingback}
\end{equation}
by the definition (\ref{sigstardef}) of $\sigma^*$. Clearly the
restriction $\sigma^*(h)>0$ is now superfluous and can be removed.

We would like to evaluate this last sum using Proposition
\ref{mvprop}. Notice that when $p$ is a prime exceeding $g$ and not
dividing the discriminant $\Delta$ of $f$, then by the definitions of
$G$ and $\sigma^*$ we have
\begin{equation*}
G(p)\sigma^*(p) = \big( 1-\frac{\sigma(p)}p \big)^{-1} \big(
\sigma(p)-\frac{\sigma(p^2)}p \big) = \sigma(p) + O(p^{-1}),
\end{equation*}
since both $\sigma(p)$ and $\sigma(p^2)$ are bounded by $g$ by
the inequality (\ref{sandor}). This implies that
\begin{equation*}
\sum_{p<x} \frac{G(p)\sigma^*(p)\log p}p = O(1) + \sum \begin{Sb}g<p<x
\\ p\dnd\Delta\end{Sb} \big( \frac{\sigma(p)\log p}p + O\big( {\log
p\over p^2} \big) \big) = \log x + O(1)
\end{equation*}
by equation (\ref{nagel}), verifying the major hypothesis
(\ref{consteq}) of Proposition \ref{mvprop} with $\kappa=1$. Also, the
remarks following the statement of Proposition \ref{mvprop} imply that
the other hypotheses (\ref{minorconverge}) and (\ref{minorprod}) are
satisfied as well, the latter with $\beta=1$.

Consequently, we may apply Proposition \ref{mvprop} to obtain
\begin{equation}
\sum_{1\le h\le \sh^\delta} \frac{G(h)\sigma^*(h)}h = \log \sh^\delta
\prod_p \big( 1 + {G(p)\sigma^*(p)\over p} + {G(p^2)\sigma^*(p^2)\over
p^2} + \dotsb \big) \big( 1-\frac1p \big) + O(1).
\label{inlight}
\end{equation}
However, each term in this product contains the telescoping series
\begin{equation*}
\begin{split}
{G(p)\sigma^*(p)\over p} + {G(p^2)\sigma^*(p^2)\over p^2} + \dotsb &=
\big( {\sigma^*(p)\over p} + {\sigma^*(p^2)\over p^2} + \dotsb \big)
G(p) \\
&= \big( {\sigma(p)\over p} - {\sigma(p^2)\over p^2} +
{\sigma(p^2)\over p^2} - {\sigma(p^3)\over p^3} + \dotsb \big) \big(
1-\frac{\sigma(p)}p \big)^{-1} \\
&= {\sigma(p)\over p} \big( 1-\frac{\sigma(p)}p \big)^{-1} = \big(
1-\frac{\sigma(p)}p \big)^{-1} - 1,
\end{split}
\end{equation*}
and thus in light of equation (\ref{inlight}), equation
(\ref{goingback}) becomes
\begin{multline}
\sum \begin{Sb}1\le h<L^\delta \\ \sigma^*(h)>0\end{Sb}
\frac{\sigma(h)}h \prod \begin{Sb}p<z \\ p^r\exdiv h\end{Sb} \big(
1-\frac{\sigma(p)}p \big)^{-1} \big( 1-{\sigma(p^{r+1}) \over
p\sigma(p^r)} \big) = \sum_{1\le h\le \sh^\delta}
\frac{G(h)\sigma^*(h)}h \\
= \log L^\delta \big( 1 + O\big( \frac1{\log L} \big) \big) \prod_p
\big( 1-\frac{\sigma(p)}p \big)^{-1} \big( 1-\frac1p \big).
\label{unga}
\end{multline}

We are now able to establish the lemma. When we insert the expressions
(\ref{bunga}) and (\ref{unga}) into the upper bound (\ref{usedsieve}),
the infinite products in (\ref{bunga}) and (\ref{unga}) cancel each
other completely, leaving the upper bound
\begin{equation*}
\sum_{d\mid P(z)} \frac{\lmp d}d \sum \begin{Sb}1\le h<L^\delta \\
\sigma^*(h)>0\end{Sb} {\sigma(dh)\over h} \le (F(s)+O(\log^{-1/3}D))
\frac{\log L^\delta}{e^\gamma\log z} \big( 1 + O\big( \frac1{\log z}
\big) \big) \big( 1 + O\big( \frac1{\log L} \big) \big).
\end{equation*}
On rearranging the various terms, writing $1/\log z$ as $s/\log D$, and
using the hypothesis that $z\le D\le\exp(\log^3L)$ to simplify the
error terms, we obtain precisely the statement (\ref{canceleq}) of the
lemma.
\end{proof}

%% file: apply.tex
We are now ready to establish Proposition \ref{grungyprop}, using the
reformulation (\ref{italics}). Let $z$ be a parameter to be specified
later subject to $L^\delta\le z\le\xi\sh^{-\delta}$. We can write the
definition (\ref{Szdef}) of $S(z)$ as
\begin{equation*}
\begin{split}
S(z) &= \sum \begin{Sb}1\le h<L^\delta \\ \sigma^*(h)>0\end{Sb} \sum
\begin{Sb}m \\ p\mid m\implies p>z\end{Sb} \3{x-L<n\le x: f(n)=mh} \\
&= \sum \begin{Sb}1\le h<L^\delta \\ \sigma^*(h)>0\end{Sb} \sum_m
\3{x-L<n\le x: f(n)=mh} \sum_{d\mid (m,P(z))} \mu(d),
\end{split}
\end{equation*}
using the characteristic property of the M\"obius function. Then, by
the upper bound sieve property (\ref{upper}), we have
\begin{equation*}
\begin{split}
S(z) &\le \sum \begin{Sb}1\le h<L^\delta \\ \sigma^*(h)>0\end{Sb}
\sum_m \3{x-L<n\le x: f(n)=mh} \sum_{d\mid (m,P(z))} \lmp d \\
&= \sum_{d\mid P(z)} \lmp d \sum \begin{Sb}1\le h<L^\delta \\
\sigma^*(h)>0\end{Sb} \3{(n,m): x-L<n\le x,\, f(n)=mh,\, d\mid m} =
\sum_{d\mid P(z)} \lmp d S_d
\end{split}
\end{equation*}
from the definition (\ref{Sddef}) of the $S_d$.

Using the expression (\ref{sdeasy}) for the $S_d$, we see that
\begin{equation}
S(z) \le L \sum_{d\mid P(z)} \frac{\lmp d}d \sum \begin{Sb}1\le
h<L^\delta \\ \sigma^*(h)>0\end{Sb} \frac{\sigma(dh)}h + O\bigg(
\sum_{d\mid P(z)} \lmp d \sum_{1\le h<L^\delta} \sigma(dh) \bigg).
\label{aftersd}
\end{equation}
The first sum can be bounded above using the inequality (\ref{Disz}).
Moreover, the $\lmp d$ have absolute value at most 1 and are supported
on integers less than $D$, which we have set equal to $z$ in order to
apply (\ref{Disz}); thus the sum in the error term is
\begin{equation*}
\ll \sum_{1\le d<z} \sum_{1\le h<L^\delta} \sigma(dh) \le
\sum_{m<zL^\delta} d(m)\sigma(m).
\end{equation*}
Since $\sigma(p)$ satisfies the asymptotic formula (\ref{consteq})
with $\kappa=1$, it follows that $d(p)\sigma(p)$ satisfies
(\ref{consteq}) with $\kappa=2$. The other hypotheses of Proposition
\ref{mvprop} are again easily verified with $\beta=2$, and thus we
may apply the upper bound (\ref{mvcor}) to obtain
\begin{equation*}
\sum_{m<zL^\delta} d(m)\sigma(m) \ll zL^\delta\log zL^\delta.
\end{equation*}

The inequality (\ref{aftersd}) now becomes
\begin{equation}
S(z) \le L \big( {2\delta\log\sh\over\log z} \big) (1+O(\log^{-1/3}z))
+ O(zL^\delta\log zL^\delta).
\label{neededlam}
\end{equation}
We want to make the main term of this upper bound small, and so we
want to choose $z$ as large as possible (without making the error term
dominant) subject to the condition $L^\delta\le z\le\xi
L^{-\delta}$. We set $z=L^{1-\delta}\log^{-2}L$, a valid choice for
sufficiently large $L$ since $\delta<1/2$. This gives
\begin{equation}
S(z) \le L\( {2\delta\over1-\delta} + O\big( {\log\log L\over\log L}
\big) \)1 (1+O(\log^{-1/3}L)) + O(L\log^{-1}L),
\label{zpicked}
\end{equation}
which is enough to establish Proposition \ref{grungyprop} and
therefore Theorems \ref{mainthm}--\ref{scantthm}.

%It is conceivable that the range of $\delta$ for which Theorems
%\ref{mainthm}--\ref{scantthm} hold could be enlarged by improving our
%estimates of the error terms. Instead of being content with an error
%of $O(1)$ in equation (\ref{bettererror}), where we counted the number
%of integers in a short arithmetic progression, we could write the
%error explicitly in terms of greatest-integer functions and then
%attempt to apply techniques for estimating exponential sums. If
%successful, this would increase the value of $z$ that we could use in
%passing from inequality (\ref{neededlam}) to inequality
%(\ref{zpicked}), thus reducing the constant $2\delta/(1-\delta)$ in
%Proposition \ref{grungyprop}, though it is not clear to what extent.

%% file: prime.tex
In this section we outline the changes to the above method needed to
establish Theorem \ref{primethm}. The ultimate object of study will
now be
\begin{equation*}
\Psi_\Lambda(f;x,y) = \sum \begin{Sb}1\le n\le x \\ p\mid f(n)\implies
p\le y\end{Sb} \Lambda(n);
\end{equation*}
if we can show that $\Psi_\Lambda(f;x,y)\gg x$, then it is easy to
deduce that $\Phi(f;x,y)\gg x/\log x$ by a simple partial summation
argument. We will use a subscripted $\Lambda$ on the notation of the
preceding sections to denote the appropriately modified quantities,
e.g.,
\begin{equation*}
C_\Lambda(f;x,y) = x-\Psi_\Lambda(f;x,y).
\end{equation*}

Theorem \ref{primethm} is a consequence of the following proposition,
which is analogous to the special case of Proposition~\ref{grungyprop}
where $L=x$:

\begin{proposition}
Let $f(t)$ be an irreducible integer-valued polynomial that is not
identically zero modulo any prime. Let $x\ge2$ be a real number and
set $\xi=\max_{n\le x}|f(n)|$, and let $\delta$ and $\ep$ be positive
real numbers such that $\delta+\ep<1/4$. Then
\begin{equation}
C_\Lambda(f;x,\xi x^{-\delta}) \le x \big(
{2\delta\over1/2-\delta-\ep} + O(\log^{-1/3}x) \big).
\label{lampropeq}
\end{equation}
\label{lamgrungyprop}
\end{proposition}

The multiplicative functions that arise in this context are
$\sigma_\Lambda(h)$, the number of roots $b$ of $f$\mod h such that
$(b,h)=1$, and $\sigma_\Lambda^*(h)$, the multiplicative function
satisfying $\sigma_\Lambda^*(p^r) = \sigma_\Lambda(p^r) -
\sigma_\Lambda(p^{r+1})/p$ for every prime power $p^r$. For example,
we have
\begin{equation}
\sigma_\Lambda(p) = \begin{cases} \sigma(p)-1, &\rmif p\mid f(0), \\
\sigma(p), &\rmif p\dnd f(0). \end{cases}
\label{porp}
\end{equation}

The outline of the proof of Proposition \ref{lamgrungyprop} is as
follows. We define
\begin{equation*}
S_\Lambda(z) = \sum \begin{Sb}1\le h<x^\delta \\
\sigma_\Lambda^*(h)>0\end{Sb} \sum \begin{Sb}n\le x \\ (n,h)=1 \\
h\mid f(n) \\ p\mid f(n)/h\implies p>z\end{Sb} \Lambda(n),
\end{equation*}
a version of the definition (\ref{Szdef}) of $S(z)$ where each term is
weighted by $\Lambda(n)$. Then, as in Section \ref{sievesec}, we have
$S_\Lambda(\xi x^{-\delta}) = C_\Lambda(f;x,\xi x^{-\delta}) +
O(x^{\delta+\ep})$, the error coming from the few terms counted by
$C_\Lambda(f;x,\xi x^{-\delta})$ that are excluded from $S_\Lambda(\xi
x^{-\delta})$ by the additional condition $(n,h)=1$ in the second
sum. Since $S_\Lambda(z)$ is again a decreasing function of $z$, it
suffices to show that $S_\Lambda(z)$ is bounded above by the
right-hand side of equation (\ref{lampropeq}) for some $2\le z\le\xi
x^{-\delta}$.

The error term in our asymptotic formula for
\begin{equation*}
(S_\Lambda)_d = \sum \begin{Sb}1\le h<x^\delta \\
\sigma_\Lambda^*(h)>0\end{Sb} \sum \begin{Sb}n\le x \\ (n,h)=1 \\
dh\mid f(n)\end{Sb} \Lambda(n) = \sum \begin{Sb}1\le h<x^\delta \\
\sigma_\Lambda^*(h)>0\end{Sb} \sum \begin{Sb}b\mod{dh} \\ (b,h)=1 \\
dh\mid f(b)\end{Sb} \sum \begin{Sb}n\le x \\ n\equiv b\mod{dh}\end{Sb}
\Lambda(n)
\end{equation*}
will now come from errors in counting the number of primes in
arithmetic progressions rather than the number of integers. If we
define
\begin{equation*}
E(t,q) = \max_{(a,q)=1} \bigg| \frac t{\phi(q)} - \sum \begin{Sb}1\le
n\le t \\ n\equiv a\mod q\end{Sb} \Lambda(n) \bigg|,
\end{equation*}
then we have
\begin{equation}
(S_\Lambda)_d = \sum \begin{Sb}1\le h<x^\delta \\
\sigma_\Lambda^*(h)>0\end{Sb} \sum \begin{Sb}b\mod{dh} \\ (b,h)=1 \\
dh\mid f(b)\end{Sb} \big( {x\over \phi(dh)} + O(E(x,dh)) \big).
\label{twosigs}
\end{equation}
Now if $d$ is squarefree, then we can write $h'=h(d,h)$ and
$d'=d/(d,h)$, so that $h'd'=hd$ and $(h',d')=1$. Then the integers $b$
such that $(b,h)=1$ and $f(b)\equiv0\mod{dh}$ are exactly those
integers such that $(b,h')=1$ and $f(b)\equiv0\mod{h'}$
and\mod{d'}. The number of such integers $b$ with $1\le b\le h'd'$ is
$\sigma_\Lambda(h')\sigma(d')$, and so the number of terms in the
inner sum of equation (\ref{twosigs}), while clearly at most
$\sigma(dh)$, is precisely
\begin{equation*}
\sigma_\Lambda(h')\sigma(d') = \sigma_\Lambda(h(d,h))\sigma\big( \frac
d{(d,h)} \big) = {\sigma_\Lambda(h(d,h))\sigma(d) \over \sigma((d,h))}.
\end{equation*}
Therefore we can derive the asymptotic formula
\begin{equation*}
(S_\Lambda)_d = x \sum \begin{Sb}1\le h<x^\delta \\
\sigma_\Lambda^*(h)>0\end{Sb} {\sigma_\Lambda(h(d,h))\sigma(d) \over
\sigma((d,h))\phi(dh)} + O\bigg( \sum_{1\le h<x^\delta} \sigma(dh)
E(x,dh) \bigg)
\end{equation*}
for any squarefree $d$, analogous to equation (\ref{sdeasy}). Since
the $\lambda_d$ are supported on squarefree integers, this leads to
the upper bound
\begin{equation*}
S_\Lambda(z) \le x \sum_{d\mid P(z)} \lmp d \sum \begin{Sb}1\le
h<x^\delta \\ \sigma_\Lambda^*(h)>0\end{Sb}
{\sigma_\Lambda(h(d,h))\sigma(d) \over \sigma((d,h))\phi(dh)} +
O\bigg( \sum_{d\mid P(z)} \lmp d \sum_{1\le h<x^\delta}
\sigma(dh)E(x,dh) \bigg),
\end{equation*}
analogous to the inequality (\ref{aftersd}).

The double sum in the main term of this inequality is similar to the
expression treated in Lemma \ref{cancellem}. Although the inner
function of $d$ and $h$ is more complicated in this case, no major
changes are needed to the method of proof of Lemma \ref{cancellem},
and we can derive the following analogous upper bound:

\begin{lemma}
\label{lamcancellem}
Let $x\ge2$ and $0<\delta<1/4$ be real numbers. For any real numbers
$z$ and $D$ satisfying $x^\delta\le z\le D\le\exp(\log^3z)$, we have
\begin{equation*}
\sum_{d\mid P(z)} \lmp d \sum \begin{Sb}1\le h<x^\delta \\
\sigma_\Lambda^*(h)>0\end{Sb} {\sigma_\Lambda(h(d,h))\sigma(d) \over
\sigma((d,h))\phi(dh)} \le \big( {2\delta\log x\over\log D} \big)
\big( {sF(s)\over2e^\gamma} + O(s\log^{-1/3}D) \big),
\end{equation*}
where $s=\log D/\log z$.
\end{lemma}

Setting $D=z$ allows us to derive an upper bound for $S_\Lambda(z)$,
analogous to equation (\ref{neededlam}), of the form
\begin{equation}
S_\Lambda(z) \le x\big( {2\delta\log x\over\log z} \big)
(1+O(\log^{-1/3}z)) + O\bigg( \sum_{1\le m<zx^\delta}
d(m)\sigma(m)E(x,m) \bigg).
\label{neededan}
\end{equation}
Since both $d(m)\ll m^{\ep/2}$ and $\sigma(m)\ll m^{\ep/2}$ for any
$\ep>0$, the latter by the observation following equation
(\ref{sandor}), we deduce that
\begin{equation*}
\sum_{1\le m<zx^\delta} d(m)\sigma(m)E(x,m) \ll (zx^\delta)^\ep
\sum_{1\le m<zx^\delta} E(x,m).
\end{equation*}
If we choose $z=x^{1/2-\delta-\ep}$ for some $\ep<1/4-\delta$ (so that
$z>x^{1/4}>x^\delta$), then we may use the Bombieri--Vinogradov
theorem to conclude that this last sum is $\ll x^{1-\ep}$ and thus
that the latter error term in the upper bound (\ref{neededan}) is $\ll
x^{1-\ep/2}$. With this choice of $z$, the upper bound
(\ref{neededan}) then becomes
\begin{equation*}
S_\Lambda(z) \le x\big( {2\delta\over1/2-\delta-\ep} +
O(\log^{-1/3}x) \big),
\end{equation*}
which establishes Proposition \ref{lamgrungyprop} and therefore
Theorem \ref{primethm}.

One can also demonstrate that a polynomial takes an abundance of
smooth values on prime arguments in short intervals, by employing
short-interval versions of the Bombieri--Vinogradov theorem. We state
the following theorem without proof, except to remark that (\ref{tda})
below uses the work of Perelli, Pintz, and Salerno
\cite{PerPinSal:BTiSI} and that (\ref{tdb}) below uses the work of
Timofeev \cite{Tim:DoAFiSI}.

\begin{theorem}
Let $h(t)$, $g$, and $k$ be as in Theorem \ref{mainthm}. Let $x\ge2$
and $0<\theta<1$ be real numbers and set $L=x^\theta$. Define
\begin{equation*}
\Phi(h;x,L,y) = \Phi(h;x,y) - \Phi(h;x-L,y).
\end{equation*}
Then when $x$ is sufficiently large, the lower bound
\begin{equation*}
\Phi(h;x,L,x^gL^{-\delta})\gg_\theta L/\log x
\end{equation*}
holds for
\begin{equation}
0<\delta<{\theta-1/2\over(2k+1)\theta} \quad \hbox{\rm if }\theta>3/5
\label{tda}
\end{equation}
and for
\begin{equation}
0<\delta<{\theta-11/20\over(2k+1)\theta} \quad \hbox{\rm if }\theta>7/12.
\label{tdb}
\end{equation}
Under the assumption of the generalized Riemann hypothesis, the
inequality $\theta>3/5$ in {\rm(\ref{tda})} may be improved to
$\theta>1/2$.
\end{theorem}

It is also clear that Proposition \ref{lamgrungyprop} implies a
smoothness result on the values a polynomial takes on a set of primes
of positive density, analogous to Theorem \ref{scantthm}.

%% file: silly.tex
In this final section we establish Theorem \ref{sillythm}. By the same
reasoning as before, it suffices to consider the case
$h(t)=f_1(t)\dotsm f_k(t)$ where the $f_i$ are distinct irreducible
polynomials of degree $g$ that are not identically zero modulo any
prime. Given a real number $x$, let $\P_i$ be the set of primes
$p\in[\frac12x^{1/k},x^{1/k}]$ such that $f_i$ has a root\mod p, and
consider the set
\begin{equation*}
\P = \{ (p_1, \dots, p_k) \in \P_1\times\dots\times\P_k : p_i\ne p_j\,
(i\ne j) \}
\end{equation*}
of $k$-tuples of distinct primes. Each element of $\P$ gives rise to a
positive integer $n\le x$ such that $h(n)$ is \osm{x^{g-1/k}} as
follows. Choose residue classes $n_i\mod{p_i}$ such that
$f_i(n_i)\equiv0\mod{p_i}$. Since the $p_i$ are distinct, we can find
a positive integer $n\le p_1\dotsm p_k$ with $n\equiv n_i\mod{p_i}$ by
the Chinese remainder theorem. Clearly each $f_i(n)\equiv0\mod{p_i}$,
and so we can write
\begin{equation*}
h(n) = f_1(n)\dotsm f_k(n) = (p_1d_1)\dotsm(p_kd_k)
\end{equation*}
for some integers $d_i$. We have $n\le p_1\dotsm p_k\le x$, and each
$d_i=f_i(n)/p_i$ is thus $\ll n^g/p_i \ll x^{g-1/k}$.  Therefore,
$h(n)$ is \osm{x^{g-1/k}}.

To determine the number of distinct values of $n$ arising in this
manner, and hence a lower bound for $\Psi(h;x,O(x^{g-1/k}))$, we need
a lower bound for the cardinality of $\P$ and an upper bound for the
number of different elements of $\P$ that could give rise to a
particular value of~$n$. If we write $\sigma_i(h)$ for the the number
of roots of $f_i$\mod h, we have $\sigma_i(p)\le g$ for any prime
$p$, and thus \def\hf{\frac12}
\begin{equation*}
\begin{split}
\#\P_i = \sum \begin{Sb}x^{1/k}/2\le p\le x^{1/k} \\ \sigma_i(p)\ge
1\end{Sb} 1 &\ge \sum_{x^{1/k}/2\le p\le x^{1/k}} \frac{\sigma_i(p)}g
\Big( \frac{\log p}p {\hf x^{1/k}\over\log\hf x^{1/k}} \Big) \\
&= \frac {x^{1/k}}{2g\log\hf x^{1/k}} \sum_{x^{1/k}/2\le p\le x^{1/k}}
\frac{\sigma_i(p)\log p}p \gg \frac {x^{1/k}}{\log x}
\end{split}
\end{equation*}
by the asymptotic formula (\ref{nagel}). Therefore the cardinality of
$\P_1\times\dots\times\P_k$ is $\gg x\log^{-k}x$; and since there are
at most $x^{1-1/k}$ $k$-tuples in $\P_1\times\dots\times\P_k$ whose
coordinates are not distinct, we see that $\#\P\gg x\log^{-k}x$.

On the other hand, if an element $(p_1,\dots,p_k)$ of $\P$ gives rise
to a particular $n$, then certainly each $p_i$ must divide
$f_i(n)$. However, each possible $p_i$ is $\gg x^{1/k}$, while each
$f_i(n)$ is $\ll x^g$; therefore there are at most $gk$ candidates for
each $p_i$ when $x$ is sufficiently large, and hence $(gk)^k\ll1$
possible elements of $\P$ that give rise to $n$. From this we conclude
that $\Psi(h;x,O(x^{g-1/k}))\gg x\log^{-k}x$, which establishes
Theorem \ref{sillythm}, aside from having $O(x^{g-1/k})$ as the
smoothness parameter instead of $x^{g-1/k}$, which we can fix by
replacing $x$ by $cx$ for a suitably small positive constant~$c$.

This technique can also demonstrate an abundance of smooth values of
polynomials of more than one variable, and in fact the range of
smoothness can be enhanced somewhat by making use of existing results
on small solutions of congruences for these polynomials.